\documentclass{article} 
\usepackage{latexsym} 
\usepackage{amsmath}
\usepackage{amssymb} 
\usepackage{amscd} 
\usepackage{amsthm}
\usepackage{graphicx}

\begin{document} 
\newtheorem{Th}{Theorem}[section]
\newtheorem{Cor}{Corollary}[section]
\newtheorem{Prop}{Proposition}[section]
\newtheorem{Lem}{Lemma}[section]
\newtheorem{Def}{Definition}[section]
\newtheorem{Ex}{Example}[section]
\newtheorem{stw}{Proposition}[section]

\theoremstyle{definition}
 \newtheorem{Rem}{Remark}[section]


\newcommand{\bet}{\begin{Th}}
\newcommand{\ent}{\stepcounter{Cor}
   \stepcounter{Prop}\stepcounter{Lem}\stepcounter{Def}
   \stepcounter{Rem}\stepcounter{Ex}\end{Th}}


\newcommand{\bec}{\begin{Cor}}
\newcommand{\enc}{\stepcounter{Th}
   \stepcounter{Prop}\stepcounter{Lem}\stepcounter{Def}
   \stepcounter{Rem}\stepcounter{Ex}\end{Cor}}
\newcommand{\bep}{\begin{Prop}}
\newcommand{\enp}{\stepcounter{Th}
   \stepcounter{Cor}\stepcounter{Lem}\stepcounter{Def}
   \stepcounter{Rem}\stepcounter{Ex}\end{Prop}}
\newcommand{\bel}{\begin{Lem}}
\newcommand{\enl}{\stepcounter{Th}
   \stepcounter{Cor}\stepcounter{Prop}\stepcounter{Def}
   \stepcounter{Rem}\stepcounter{Ex}\end{Lem}}
\newcommand{\bef}{\begin{Def}}
\newcommand{\enf}{\stepcounter{Th}
   \stepcounter{Cor}\stepcounter{Prop}\stepcounter{Lem}
   \stepcounter{Rem}\stepcounter{Ex}\end{Def}}
\newcommand{\ber}{\begin{Rem}}
\newcommand{\enr}{
   \stepcounter{Th}\stepcounter{Cor}\stepcounter{Prop}
   \stepcounter{Lem}\stepcounter{Def}\stepcounter{Ex}\end{Rem}}
\newcommand{\bee}{\begin{Ex}}
\newcommand{\ene}{
   \stepcounter{Th}\stepcounter{Cor}\stepcounter{Prop}
   \stepcounter{Lem}\stepcounter{Def}\stepcounter{Rem}\end{Ex}}
\newcommand{\Proof}{\noindent{\it Proof}\/: }

\newcommand{\EE}{\mbox{\bf{E}}}
\newcommand{\QQ}{\mbox{\bf{Q}}}
\newcommand{\R}{\mbox{\bf{R}}}
\newcommand{\C}{\mbox{\bf{C}}}
\newcommand{\ZZ}{\mbox{\bf{Z}}}
\newcommand{\NN}{\mbox{\bf{N}}}
\newcommand{\PP}{\mbox{\bf{P}}}
\newcommand{\uuu}{\boldsymbol{u}}
\newcommand{\xxx}{\boldsymbol{x}}
\newcommand{\aaa}{\boldsymbol{a}}
\newcommand{\bbb}{\boldsymbol{b}}
\newcommand{\AAA}{\mbox{\bf{A}}}
\newcommand{\BBB}{\mbox{\bf{B}}}
\newcommand{\ccc}{\boldsymbol{c}}
\newcommand{\iii}{\boldsymbol{i}}
\newcommand{\jjj}{\boldsymbol{j}}
\newcommand{\kkk}{\boldsymbol{k}}
\newcommand{\rrr}{\boldsymbol{r}}
\newcommand{\FFF}{\boldsymbol{F}}
\newcommand{\yyy}{\boldsymbol{y}}
\newcommand{\ppp}{\boldsymbol{p}}
\newcommand{\qqq}{\boldsymbol{q}}
\newcommand{\nnn}{\boldsymbol{n}}
\newcommand{\vvv}{\boldsymbol{v}}
\newcommand{\eee}{\boldsymbol{e}}
\newcommand{\fff}{\boldsymbol{f}}
\newcommand{\www}{\boldsymbol{w}}
\newcommand{\0}{\boldsymbol{0}}
\newcommand{\lon}{\longrightarrow}
\newcommand{\ga}{\gamma}
\newcommand{\pa}{\partial}
\newcommand{\QED}{\hfill $\Box$}
\newcommand{\id}{{\mbox {\rm id}}}
\newcommand{\Ker}{{\mbox {\rm Ker}}}
\newcommand{\grad}{{\mbox {\rm grad}}}
\newcommand{\ind}{{\mbox {\rm ind}}}
\newcommand{\rot}{{\mbox {\rm rot}}}
\newcommand{\diver}{{\mbox {\rm div}}}
\newcommand{\Gr}{{\mbox {\rm Gr}}}
\newcommand{\Diff}{{\mbox {\rm Diff}}}
\newcommand{\Symp}{{\mbox {\rm Symp}}}
\newcommand{\Ct}{{\mbox {\rm Ct}}}
\newcommand{\Uns}{{\mbox {\rm Uns}}}
\newcommand{\rank}{{\mbox {\rm rank}}}
\newcommand{\sign}{{\mbox {\rm sign}}}
\newcommand{\supp}{{\mbox {\rm supp}}}
\newcommand{\Spin}{{\mbox {\rm Spin}}}
\newcommand{\Sp}{{\mbox {\rm sp-codim}}}
\newcommand{\Int}{{\mbox {\rm Int}}}
\newcommand{\Hom}{{\mbox {\rm Hom}}}
\newcommand{\codim}{{\mbox {\rm codim}}}
\newcommand{\ord}{{\mbox {\rm ord}}}
\newcommand{\Iso}{{\mbox {\rm Iso}}}
\newcommand{\corank}{{\mbox {\rm corank}}}
\def\mod{{\mbox {\rm mod}}}
\newcommand{\pt}{{\mbox {\rm pt}}}
\newcommand{\enP}{\hfill $\Box$ \par\vspace{5truemm}}
\newcommand{\spe}{\vspace{0.4truecm}}
\newcommand{\im}{\mathrm{i}}

\title{Classification of phase singularities \\ 
for complex scalar waves} 

\author{
Jiro ADACHI\thanks{\scriptsize 
Partially supported by 
Grants-in-Aid for Young Scientists (B), No.~17740027, 
21st Century COE Program ``Topological Science and Technology''
Hokkaido University, 
and 21st Century COE Program 
``Mathematics of Nonlinear Structures via Singularities'' Hokkaido University. 
} 
\ and \ 
Go-o ISHIKAWA\thanks{\scriptsize 
Partially supported by 
Grants-in-Aid for Scientific Research, No.~14340020, 
21st Century COE Program ``Mathematics of Nonlinear Structures via Singularities'', 
and 21st Century COE Program ``Topological Science and Technology'' 
Hokkaido University.
}
}

\renewcommand{\thefootnote}{\fnsymbol{footnote}}
\footnotetext{\scriptsize 
Key words: wave dislocation, equi-phase portrait, 
optical vortex, Helmholtz equation.}
\footnotetext{\scriptsize
2000 {\it Mathematics Subject Classification}\/:  
Primary 58K40, 78A40; Secondly 78A05. }

\date{}
\maketitle

\begin{abstract}
Motivated by the importance and universal character of phase  
singularities which are clarified recently,
we study the local structure of equi-phase loci near the
dislocation locus of complex valued planar and spatial waves, from
the viewpoint of singularity theory of differentiable mappings,
initiated by H.~Whitney and R.~Thom.
The classification of phase-singularities
are reduced to the classification of planar curves by radial  
transformations
due to the theory of A.~du Plessis, T.~Gaffney and L.~Wilson.
Then fold singularities are classified into hyperbolic and elliptic  
singularities.
We show that the elliptic singularities are never realized by any  
Helmholtz waves, while the hyperbolic singularities are realized in fact.
Moreover, the classification and realizability of Whitney's cusp, 
as well as its bifurcation problem are considered in order to explain
the three points bifurcation of phase singularities.
In this paper, we treat the dislocation of linear waves mainly,  
developing the basic and universal method, 
the method of jets and transversality,
which is applicable also to non-linear waves.
\end{abstract}

\section{Introduction.}

A complex scalar wave has the locus, the dislocation locus, 
where its phase is not defined. The local 
structure of equi-phase loci near the dislocation locus is called a 
{\it phase-singularity} \cite{Nye2}. 
The phase singularities are called {\it optical vortices}\/ in optics 
and are very basic and important objects in any science related to waves and 
quanta. 
In this paper we give the exhaustive classification of phase 
singularities of complex scalar waves of low codimension. 

In \cite{Nye1}\cite{Nye2}, J.~F.~Nye 
constructed extensively complex scalar global 
planar waves satisfying the Helmholtz equation, 
with detailed analysis of those examples. Also he gave, by his examples,  
an explanation of an experimental bifurcation process of phase singularities: 
one degenerate singular point bifurcates to three singular points 
and then another singular point annihilates with one of the three.  
He intends to explore the phase singularities from 
an analogy with the catastrophe theory \cite{Thom}\cite{Arnold}. 

In this paper, we understand phase singularities clearly 
using the singularity theory of differentiable mappings 
\cite{AGV}\cite{Mather3}\cite{Wall}\cite{Damon}. 

The planar complex scalar wave  
can be regarded, from the general point of view, simply as 
 a differentiable mapping from the 
plane to the plane of complex numbers. 
Then, by a theorem of H.~Whitney \cite{Whitney}, 
the generic singularities of the wave, as a differentiable mapping, 
are just the {\it fold singularities}\/ and the {\it cusp singularities}. 
The singular values form on the plane of complex numbers 
an immersed curve, the discriminant, with several number of cusps. 
Then generically the discriminant does not hit the zero, 
so that the zero is a regular value; generic phase singularities 
are regular, namely, 
locally diffeomorphic to the standard radial lines emitted from the origin. 
However, for a generic time-depending wave, the curve of 
singular values moves and momentarily may hit the zero. 
Thus, generically momentary wave can have degenerate phase 
singularities described by the fold singularities. 
Moreover, for a generic two parameter family of plane wave, 
the cusp singularity occurs as more complicated phase singularities. 

The above simplified story must be examined twofold: 
First, in Whitney's theorem, the singularities are classified 
by means of arbitrary local diffeomorphisms of the source plane 
and the target plane. 
However, for the classification of phase singularities, 
we concern with the equi-phase lines and thus 
need to consider finer classification using only diffeomorphisms 
which preserve the radial lines on the target plane. 
Second, because waves must obey several natural conditions 
given by, say, the Helmholtz equations and the wave equations,  
more than just the differentiability, we must consider the realizability 
of singularities and determine generic 
singularities among waves satisfying those conditions. 

We clarify the equivalence relations for phase singularities, and 
thus classify all phase singularities of low codimension, 
and discuss the realizability by the Helmholtz waves of phase singularities. 
Further, 
we propose the new explanation for the experimental bifurcation process 
treated in \cite{Nye1}\cite{Nye2}.

\ 

In the next section, we formulate our equivalence relation 
providing the base of our classification. 
A natural and refined classification by radial transformations 
is established on phase singularities 
for planar and spatial complex scalar waves. 
Then, we give the exact classification of generic complex planar waves 
without conditions motivated from physics. 

In \S 3, the realization of singularities by the Helmholtz waves is examined 
by concrete examples, which have a different character 
with Nye's examples in \cite{Nye1}. 

In \S 4, we treat phase singularities of 
spatial complex scalar waves and 
consider their realizability. 

In \S\ref{sec:planar_curves}, 
the classification problem of planar waves is reduced 
to that of planar curves under diffeomorphisms preserving 
radial lines. 

We introduce in \S\ref{sec:HelmholtzJet} the notion of Helmholtz jet spaces 
and transversality to discuss genericity of singularities for Helmholtz waves.

In \S\ref{sec:Schrodinger}, 
as an application of the method developed in this paper, we discuss 
the bifurcation problem of phase singularities of solutions to 
non-linear Schr{\" o}dinger equations. 


In this paper, we consider local classification problem of phase singularities.
For the global topology of dislocation locus, see 
\cite{BD1}\cite{BD2}. 

For other applications of the singularity theory to 
solutions of partial differential equations, 
see \cite{Iz}\cite{IM}\cite{Damon2} for instance.

\section{Phase singularities for planar complex scalar waves.} 
\label{sec:PhaseSingPlanar}

We denote by $\C$ the plane of complex numbers and write a complex number 
as $u + \im w = re^{\im \theta}$, 
$u, w$ being the real part and the imaginary part respectively, 
while $r, \theta$ the modulus (or the amplitude) 
and the argument (or the phase) respectively. 

Let us consider a complex scalar wave 
$$
\Psi = \Psi(x, y, t) 
= u(x, y, t) + \im w(x, y, t)
$$ 
on the $(x,y)$-plane 
depending on the time (or any other one-parameter). 
First we regard $\Psi$ as just a time-depending complex valued function. 
We assume $u(x, y, t)$, $v(x, y, t)$ 
are differentiable (i.e. $C^\infty$) functions.

If $\Psi(x, y, t) \not= 0$ at a point $(x, y)$ and at a moment $t$, 
then we can 
write $\Psi(x, y, t) = r(x, y, t)e^{\im \theta(x, y, t)}$ uniquely 
with $r(x, y, t) > 0$ and $\theta(x, y, t)$ mod. $2\pi$. 
Then, we are concerned with the wave dislocation locus at a moment $t = t_0$ 
$$
\{(x, y) \mid \Psi(x, y, t_0) = 0\} = \{(x, y) \mid 
 u(x, y, t_0) = 0, w(x, y, t_0) = 0\}
$$ 
and the equi-phase curves $\{(x,y)\mid\theta(x, y, t_0) = \text{const.}\}$ outside of 
the wave dislocation locus. 

Then, in the framework of singularity theory of differentiable mappings, 
we introduce 
the notion of {\it radial transformations}\/ and 
give the exact classification results of singularities relatively to the 
radial transformations. 

A {\it radial transformation}\/ on $\C$ near $0$ 
is a diffeomorphism, an invertible differentiable transformation, 
$\tau(u, w) = (U, W)$, $\tau\colon (\C, 0) \to (\C, 0)$ 
which sends any radial line $\{ \theta = {\rm const.}\}$ to a radial line. 
In fact, a diffeomorphism $\tau(u, w) = (U, W)$ is a radial transformation 
if and only if 
there exists a positive function $\rho(u, w)$ and real numbers $a, b, c, d$ 
with $ad - bc \not= 0$ such that 
$$
U = \rho(u, w)(au + bw), \qquad W = \rho(u, w)(cu + dw). 
$$

For the classification, we define the equivalence relation on phase singularities: 
Two functions $\Psi(x, y, t)$ and 
$\Phi(x, y, t) = u'(x, y, t) + \im w'(x, y, t)$ 
are {\it radially equivalent}\/ at points and moments $(x_0, y_0, t_0)$ and 
$(x_0', y_0', t'_0)$ respectively, 
if there exist a local diffeomorphism $\sigma(x, y) = (X(x, y), Y(x, y))$ on the plane 
with $X(x_0, y_0) =  x_0', Y(x_0, y_0) = y_0'$ and 
a local radial transformation $\tau(u, w) = (U(u, w), W(u, w))$ near the origin on $\C$ 
such that 
$$
\begin{array}{rcl}
\vspace{0.2truecm} 
u(X(x, y), Y(x, y), t_0) & = & U(u'(x, y, t'_0), w'(x, y, t'_0)), \\ 
w(X(x, y), Y(x, y), t_0) & = & W(u'(x, y, t'_0), w'(x, y, t'_0)), 
\end{array}
$$
namely that $\Psi(\sigma(x, y), t_0) = \tau(\Phi(x, y, t'_0))$. 

\begin{figure}[htb] 
 \centering
 \includegraphics[width=10cm,keepaspectratio,clip]{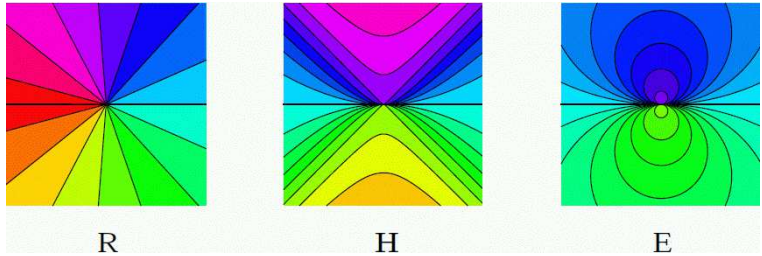}
 \caption{phase singularities}
 \label{fig:phasesings}
\end{figure}

\bet \label{classif2d}
For a generic complex valued function $\Psi(x, y, t)$, 
the phase singularity at any point and any moment 
$(x_0, y_0, t_0)$ is equivalent under radial transformations to 
the regular singularity 
$$
{\mbox{\rm R:}} \quad \psi(x, y) = x  + \im y,  
$$
the hyperbolic singularity
$$
{\mbox{\rm H:}} \quad \psi(x, y) = x^2  - y^2  + \im y,  
$$
or to the elliptic singularity
$$
{\mbox{\rm E:}} \quad \psi(x, y) = x^2  + y^2 + \im y,  
$$
at the origin $(x, y) = (0, 0)$. 
\textup{(}see \textup{Figure}~$\ref{fig:phasesings}$. \textup{)}
\ent

Each phase singularity of the classification in Theorem~\ref{classif2d} is  
determined
by its two jet actually.

\ber
Besides the wave dislocation, we can classify generic critical points of 
phase functions defined outside of the dislocation locus: 
The generic critical points 
are the {\it non-degenerate maximal(minimal) points}, 
the {\it saddle points}\/ and the {\it cuspidal points}. 
The last one bifurcates to 
one maximal(minimal) point and one saddle point.  
\enr

Both hyperbolic and elliptic singularities are equivalent 
to the {\it fold singularity}
$$
\psi : (x, y) \mapsto (u, w) = (x^2, y)
$$
under arbitrary diffeomorphisms not necessarily radial transformation, 
namely, under the right-left equivalence. 

For a momentary complex wave 
$\psi(x, y) = u(x, y) + \im w(x, y)$ on the plane, 
the locus in $\C$ of complex values $\psi(x_0, y_0)$ for $(x_0, y_0)$ with 
$$
\det
\left(
\begin{array}{cc}
\dfrac{\pa u}{\pa x} & \dfrac{\pa u}{\pa y} \\
\dfrac{\pa w}{\pa x} & \dfrac{\pa w}{\pa y} 
\end{array}
\right) (x_0, y_0) = 0, 
$$
is called the {\it discriminant\/} of the complex wave $\psi$. 

Note that the above normal form of fold singularity 
is never generic as the phase singularity. 
In fact, the discriminant of $\psi$ in that case 
is the $w$-axis in $\C$ which has the infinite tangency (actually coincides) 
with the radial lines $\{ \theta = \pi/2 \}$ and $\{ \theta = 3\pi/2 \}$. 
Generically the discriminant must be tangent to the radial lines 
in non-degenerate manner, namely, in the second order tangency. 
Then there are two possibility of non-degenerate tangency of the discriminant of fold singularities at $0 \in \C$; 
the image of $\psi$ (the value set of $\psi$) 
is concave or convex. These correspond, respectively, 
to the hyperbolic singularity 
and elliptic singularity. 

\begin{figure}[htb] 
 \centering
 \includegraphics[width=10cm,keepaspectratio,clip]{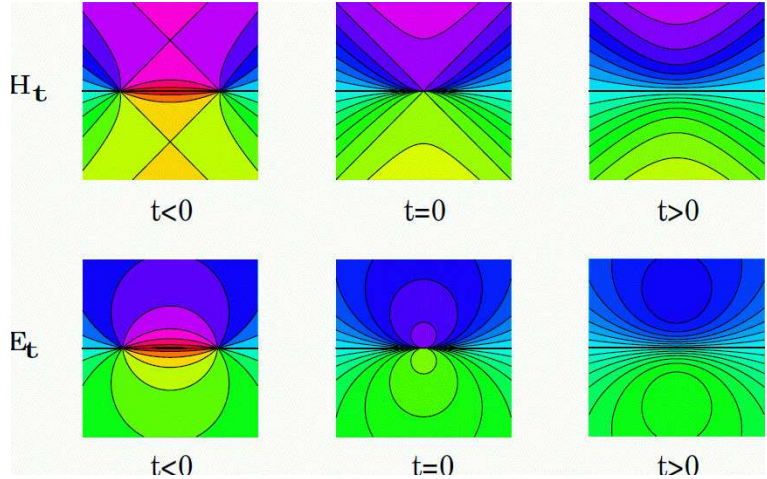}
 \caption{The bifurcations of the hyperbolic phase singularity (top) and 
the elliptic phase singularity (bottom).}
 \label{fig:bifurcations}
\end{figure}

Moreover, the generic bifurcations on $t$ of the hyperbolic singularities  
and the elliptic singularities are given by 
$$
\begin{array}{rrcl}
\vspace{0.2truecm}
{\mbox{{\rm H}$_t$}} : & 
\Psi(x, y, t) & = & x^2  - y^2 + t + \im y, \quad (t \in \R), \\  
{\mbox{{\rm E}$_t$}} : & 
\Psi(x, y, t) & = & x^2  + y^2 + t +  \im y, \quad (t \in \R). 
\end{array}
$$
(see Figure~$\ref{fig:bifurcations}$.)

The picture of the bifurcation of the hyperbolic singularity 
can be seen in Fig.~6 of \cite{DV}.

\ber
The classification of phase singularities is closely related to 
the classification of plane curves under diffeomorphisms preserving 
a given singular foliation on the plane. 
Then, one of the most delicate cases is the case 
when the foliation is formed by radial lines, which is given by 
Euler vector field $X = x\dfrac{\pa}{\pa x} + y\dfrac{\pa}{\pa y}$ 
(cf. \cite{Zhitomirskii}). 
That is the case we are treating in this paper (see \S 5). 
\enr

The discriminant of the fold singularity is a regular curve. 
Then degenerate phase singularities are classified as follows: 

\bep
The phase singularities arising from fold singularities are classified into 
$$
\psi_m(x, y) = x^2 \pm y^m + \im y, \quad m = 2, 3, 4, \dots, 
$$
under radial transformations {\rm (}and diffeomorphisms on the target{\rm )}, 
provided the discriminant curve has a contact 
with the tangent line at the origin in a finite multiplicity. 
\enp

\ber
The generic bifurcation of $\psi_m(x, y)$ 
is described by the family 
$$
x^2 \pm y^m + t_{m-1}y^{m-1} + t_{m-2}y^{m-2} +  \cdots 
+ t_2 y^2 + t_0 + \im y,
$$
with $(m - 1)$-parameters $t_0, t_2, \dots, t_{m-1}$. 
\enr

A momentary complex wave $\psi\colon \left(\R^2, (x_0, y_0)\right) \to (\C, 0)$ 
is called a {\it Whitney's cusp}\/ or simply a {\it cusp\/} 
if it is right-left equivalent (under local diffeomorphisms on 
$\R^2$ and $\C$ which are not necessarily radial) to the mapping 
$\psi(x, y) = x^3 + xy + \im y$. 
We are interested in this type of phase singularity 
because there occurs a three points bifurcation 
by just a translation $\psi_a(x, y) = x^3 + xy + \im (y + a), (a \in \R)$. 
(For the classification of more degenerate singularities 
under the right-left equivalence relations, 
see \cite{Rieger}\cite{RR}). 

The Whitney's cusp appears generically in two parameter families 
of planar complex valued functions. 

Then, we have 
\bep
{\rm (The radial classification of Whitney's cusps):} 
Any Whitney's cusp is equivalent under radial transformations to the standard 
function  
$\psi(x, y) = x^3 + xy + \im y$. 
\enp

\smallskip

\begin{figure}[htb] 
 \centering
 \includegraphics[width=10cm,keepaspectratio,clip]{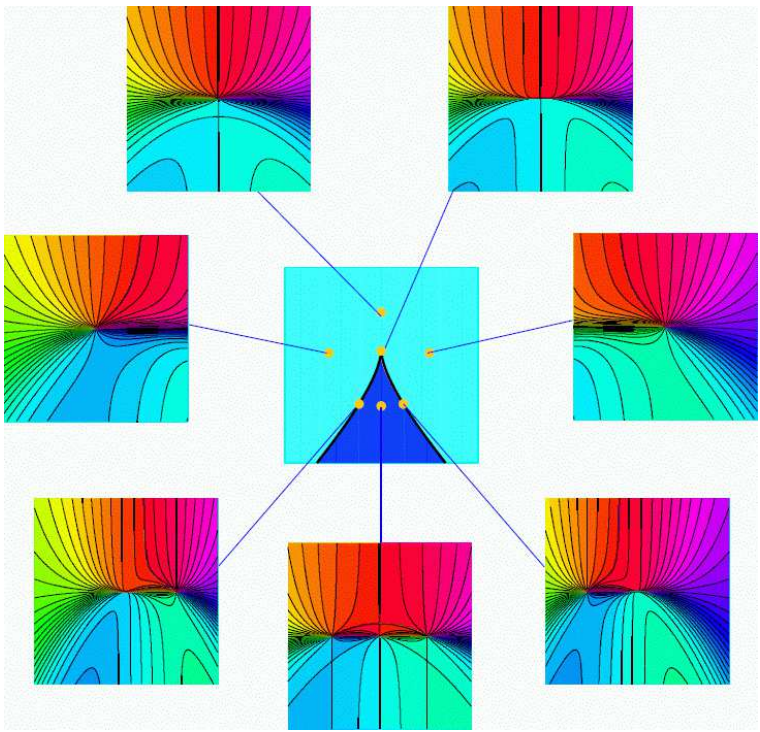}
 \caption{The two parameter bifurcation of the cusp phase singularity. }
 \label{fig:cuspbifu}
\end{figure}

The typical bifurcation of the phase singularities for a Whitney's cusp 
is described by 
$\psi_{a, b}(x, y) = x^3 + xy + b + \im (y + a), \ (a, b \in \R)$. 
\textup{(}see \textup{Figure}~$\ref{fig:cuspbifu}$\textup{)}

\ber
\textrm{
The bifurcation problem of phase singularities arising from Whitney cusps 
is related to web geometry (\cite{AkGol}, \cite{Dufour}). 
 In fact, generic two parameter families of Whitney cusps 
define $3$-webs on the plane, 
and their classification by radial transformations provides functional moduli
(Remark~\ref{3-web})}. 
\enr

\section{Phase singularities of the Helmholtz waves.}

Now, we ask the physical reality of the classification; 
the instantaneous appearance of singularities for wave functions 
satisfying the wave equation and the Helmholtz equation. 
Namely, we assume the wave $\Psi(x, y, t)$ satisfies the 
wave equation 
$$
\dfrac{\pa^2 \Psi}{\pa t^2} = c^2\nabla^2\Psi,  
$$
for a positive real number $c$, 
where $\nabla^2 = \dfrac{\pa^2}{\pa x^2} + \dfrac{\pa^2}{\pa y^2}$ 
is the Laplacian. 
Moreover we assume
that $\Psi(x, y, t)$ satisfies the Helmholtz equation 
$$
\nabla^2\Psi + k^2\Psi = 0, 
$$ 
for a positive real number $k$ as a very natural physical assumption 
for monochromatic waves.  
We call a function which satisfies the Helmholtz equation 
the \emph{Helmholtz function}. 
Further, we call the Helmholtz function with a parameter $t$ 
which satisfies the wave equation as well the \emph{Helmholtz wave}. 
Note that if we have a solution $\Psi$ for $c = 1$, $k = 1$, 
then, by setting 
$\widetilde{\Psi}(x, y, t) = \Psi(k x, k y, \frac{c}{k} t)$,  
we have a solution $\widetilde{\Psi}$ for general $c$ and $k$. 

By solving the Cauchy problem properly, we obtain the following. 

\bep
The complex valued function 
$$
\psi(x, y) = \cos y - \cos x  + \im \sin y
$$
satisfies Helmholtz equation $\psi_{xx} + \psi_{yy} + \psi = 0$ and 
has the hyperbolic singularity at the origin. 
Moreover the hyperbolic singularity with its generic bifurcation 
is realized by a Helmholtz wave 
$$
\Psi(x, y, t) = (\cos y - \cos x  + \im \sin y)\cos t + \cos y \sin t, 
$$
{\rm (}for $k = 1, c = 1${\rm )}. 
\enp

We see it is radially equivalent to the normal form
simply by observing its Taylor expansion.

To the contrary, we observe:  

\bep \label{prop:realize_ellip}
Elliptic singularities are not realized 
as a function satisfying the Helmholtz equation. 
\enp

\Proof 
Suppose a function $\psi(x, y)$ is radially equivalent 
to the elliptic singularity. 
Then, the image of $\psi$ is convex at the origin 
where the tangent line supports. 
Suppose the function $\psi(x,y)$ satisfies the Helmholtz equation with $k=1$. 
From the equation $\psi_{xx} + \psi_{yy} + \psi = 0$, 
we see the Hessian of $\psi$ is traceless 
at the dislocation locus $\{ \psi = 0 \}$, 
so are the real part ${\rm{Re}}({\rm{Hess}}\psi)$ 
and the imaginary part ${\rm{Im}}({\rm{Hess}}\psi)$. 
Thus, for any real numbers $\lambda, \mu$, 
$$
\lambda{\rm{Re}}({\rm{Hess}}\psi) + \mu{\rm{Im}}({\rm{Hess}}\psi)
$$
is never a definite matrix. 
However, the linear projection along the tangent line to the image of $\psi$ 
must be definite. 
This leads to a contradiction.  
\qed

\

In fact, as for the generic classification of phase singularities 
for one-parameter families of complex valued functions 
satisfying Helmholtz equation, 
 we have: 

\bet
\label{thm:Helmholtz_generic_fct}
The generic phase singularities of planar Helmholtz functions are 
regular singularities and hyperbolic singularities. 
\ent

\bet
\label{Helmholtz generic}
The generic phase singularities of planar Helmholtz waves are 
regular singularities and hyperbolic singularities. 
\ent
\noindent
Theorem~\ref{thm:Helmholtz_generic_fct}, and Theorem~\ref{Helmholtz generic} 
are proved in Section~\ref{sec:HelmholtzJet}. 

For the cusp singularities, we have: 

\bep
\label{real cusp}
A Whitney's cusp is realized as a Helmholtz wave. 
In fact, 
$$
\psi(x, y) = x^3\cos y + (x - 3 xy)\sin y + \im \sin y
$$ 
is a Whitney's cusp 
satisfying Helmholtz equation {\rm (}$k = 1${\rm )}\textup{:} 
$\psi_{xx} + \psi_{yy} + \psi = 0$. 
Moreover, 
$$
\Psi(x, y, t) = (x^3\cos y + (x - 3 xy)\sin y + \im \sin y)\cos t 
+ \im\cos y\sin t
$$
gives a deformation of $\psi$ by a Helmholtz wave 
\textup{(}$k = 1, c = 1$\textup{)} 
describing a three point bifurcation of the phase singularity. 
\enp

\ber
By a similar construction to Proposition~\ref{real cusp}, 
we have another realization 
$$
\psi(x, y) = x^2 \cos y - y\sin y + \im \sin y
$$
of hyperbolic singularities. 
\enr

\ber 
Apart from the classification problem of phase singularities, 
we can show that 
generic Helmholtz function $\psi : \R^2 \to \C$ is, locally at any point 
in $\R^2$, {\it right-left equivalent}\/ 
to a regular point, to a fold point or to a cusp point. 
\enr

\section{The radial classification of 
phase singularities for spatial complex scalar waves.} 

We study, in this section, the phase singularities of spatial waves 
$\Psi = \Psi(x, y, z, t) : \R^3\times \R \to \C$. 

The generic singularities of differentiable mappings 
$\R^3 \to \C$ consist of the {\it definite fold singularities}, 
the {\it indefinite fold singularities}\/ and the {\it cusp singularities}. 
The normal forms of them are given by 
$$
\begin{array}{rccl}
{\rm the \ definite \ fold \ singularity :} & 
\psi(x, y, z) & = & x^2 + y^2 + \im z, \\
{\rm the \ indefinite \ fold \ singularity :} & 
\psi(x, y, z) & = & x^2 - y^2 + \im z, \\
{\rm the \ cusp \ singularity :} & 
\psi(x, y, z) & = & x^3 + xy + z^2 + \im y, \\
\end{array}
$$
under the left-right equivalence \cite{GG}. 
For a generic complex valued function $\Psi(x, y, z, t)$, 
only fold singularities may appear as a phase singularity. 
Moreover the discriminant curve 
has non-degenerate tangency with the tangent line at the origin. 
Thus we have 

\begin{figure}[htb] 
  \centering
  \includegraphics[width=10cm]{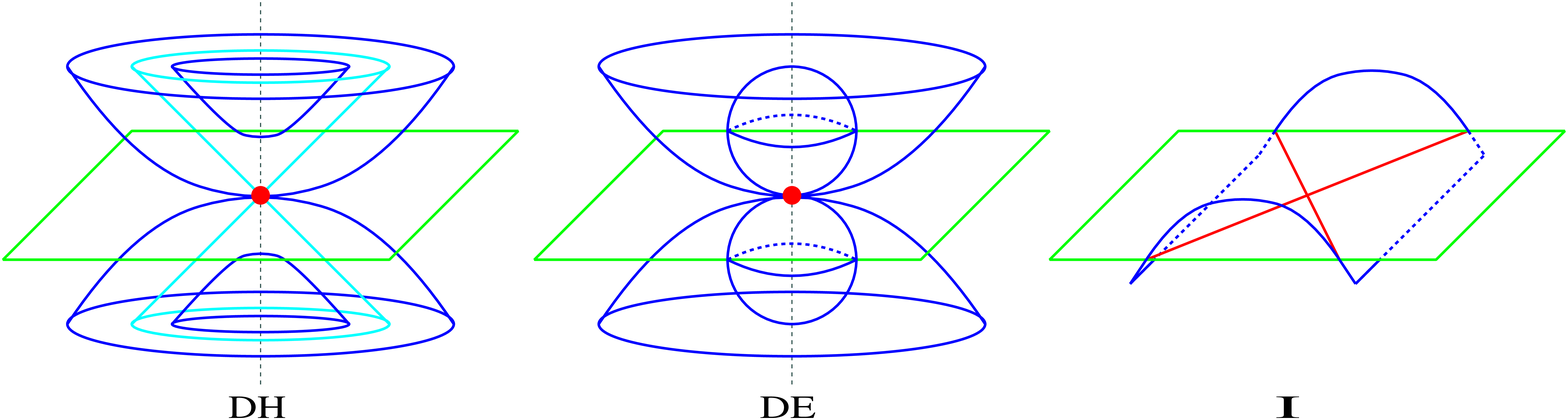}  
  \caption{Phase singularities for spatial waves}
  \label{fig:3Dphasesing}
\end{figure}

\bet
For a generic spatial complex valued function $\Psi(x, y, z, t)$, 
the phase singularity at any point and any moment 
$(x_0, y_0, z_0, t_0)$ is equivalent, under the radial transformation, to 
the regular singularity 
$$
{\mbox{\rm R:}} \quad \psi(x, y, z) = x  + \im y,  
$$
to the definite hyperbolic singularity
$$
{\mbox{\rm DH:}} \quad \psi(x, y, z) = x^2 + y^2 - z^2 + \im z,  
$$
to the definite elliptic singularity
$$
{\mbox{\rm DE:}} \quad \psi(x, y, z) = x^2 + y^2 + z^2 + \im z ,  
$$
or to the indefinite singularity
$$
{\mbox{\rm I:}} \quad \psi(x, y, z) = x^2  - y^2 - z^2 + \im z,  
$$
at the origin $(x, y, z) = (0, 0, 0)$. 
\textup{(}See \textup{Figure~\ref{fig:3Dphasesing}}.\textup{)}
\ent

\begin{figure}[htb] 
  \centering
  \includegraphics[width=10cm]{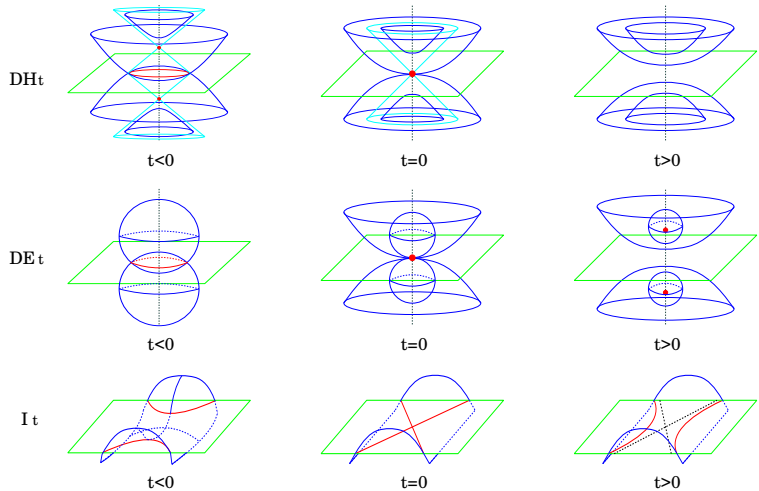}
  \caption{The bifurcations of phase singularities of spatial scalar waves. }
  \label{fig:spbif}
\end{figure}

  Similarly to the case of planar complex scalar waves, 
the generic bifurcations on $t$ of the definite hyperbolic singularities, 
the definite elliptic singularities, and the indefinite singularities 
are given by
\begin{eqnarray*}
  &\mathrm{DH_t}&:\ \Psi(x,y,z,t)=x^2+y^2-z^2+t+\im z, \\
  &\mathrm{DE_t}&:\ \Psi(x,y,z,t)=x^2+y^2+z^2+t+\im z, \\
  &\mathrm{I_t}&:\ \Psi(x,y,z,t)=x^2-y^2-z^2+t+\im z. 
\end{eqnarray*}
(see Figure~\ref{fig:spbif}.)

For cusp singularities, we have 

\bep
\label{Whitney cusp}
The phase singularities which come from cusp singularities are 
all radially equivalent to 
$$
\psi(x, y, z)  =  x^3 + xy + z^2 + \im y. 
$$
\enp 
\noindent Proposition~\ref{Whitney cusp} is induced 
from Lemma~\ref{simple cusp}. 

A complex valued function $\Psi(x,y,z,t)$ is called a \emph{Helmholtz wave}\/ 
if it satisfies the wave equation 
$\Psi_{tt}=c^2(\Psi_{xx}+\Psi_{yy}+\Psi_{zz})$ for a positive real number $c$ 
and the Helmholtz equation $\Psi_{xx}+\Psi_{yy}+\Psi_{zz}+k^2\Psi=0$ 
for a positive real number $k$. 

As for the realizability of the spatial waves as Helmholtz waves, we have: 

\bep
The definite hyperbolic singularity and the indefinite singularity 
together with their generic bifurcations are realized by 
Helmholtz waves {\rm (}for $k = 1, c = 1${\rm )}: 
$$
\begin{array}{ccc}
\vspace{0.5truecm}
{\rm DH_t} & {\rm : } & \Psi(x, y, z, t) = 
(- \cos x - \cos y + 2 \cos z  + \im \sin z)\cos t + \cos z \sin t. \\
{\rm I_t} & {\rm :} &  \Psi(x, y, z, t) = (- 2\cos x + \cos y +  \cos z + \im \sin z)\cos t 
+ \cos z \sin t. 
\end{array}
$$
\enp

Moreover, in a similar way to Proposition~\ref{prop:realize_ellip}, we have 
\bep
Definite elliptic singularities are not realized as Helmholtz waves. 
\enp

The cusp singularity is realized as a Helmholtz function: 
\bep
The complex valued function 
$$
\psi(x, y, z) = 
x^3\cos y + (x - 3xy)\sin y - \cos y + \cos z + \im\sin y
$$ 
satisfies Helmholtz equation $\psi_{xx} + \psi_{yy} + \psi_{zz} 
+ \psi = 0$ and is radially equivalent to the cusp singularities. 
\enp

\section{Radial classification of planar curves.} \label{sec:planar_curves}

The classification problem of complex waves under radial transformations 
is reduced to the classification problem 
of planar curves under radial transformations, 
by means of du Plessis, Gaffney, Wilson's theory \cite{dPGW}\cite{GW}. 
The theory reduces the classification to that of discriminants 
with an exception. 
The exceptional cases are hyperbolic and elliptic singularities. 
Although they have the same discriminants, they are not radially equivalent. 
It depends on whether the image is convex or concave. 
The following results on curves have been applied 
to the classification of discriminants, and then, phase singularities. 

\bel
{\rm (The radial classification of regular curves): }
A regular curve through the origin on $\C$ is transformed 
by radial transformations 
to the curve $u = w^m$ for some integer $m \geq 2$ or 
$u = g(w)$ for some function with null derivatives 
$g^{(i)}(0) = 0, i = 0, 1, 2, 3, \dots$. 
\enl

\Proof
Using a linear transformation, we may suppose the regular curve is given by 
$u = g(w)$ for a function $g(w)$. Suppose $\ord g = 0$ at $w = 0$ and 
$g^{(m)}(0) > 0$. 
Then we can write $u = (a(w)w)^m$ for a function $a(w)$ with $a(0) > 0$. 
Set $\rho(w) = a(w)^{\frac{m}{m-1}}$, and define the radial transformation 
$U = \rho(w)u,\ W = \rho(w)w$. 
Then $U = \rho(w)\{a(w)W/\rho(w)\}^m = W^m$. 
\qed

In Section~\ref{sec:PhaseSingPlanar}, 
we study deformations of phase singularities. 
For them we observe the following. 
\ber
Let $(u(t, \lambda), w(t, \lambda))$ be a deformation of 
the curve $(u(t, 0), w(t, 0)) = (t^m, t)$: 
$$
\begin{array}{rcl}
u(t, \lambda) & = & \alpha_0(\lambda) + \alpha_1(\lambda)t + 
\cdots + \alpha_m(\lambda)t^m + \cdots, \\
w(t, \lambda) & = & 
\beta_0(\lambda) + \beta_1(\lambda)t + \cdots, 
\end{array}
$$ 
with $\alpha_0(0) = \cdots = \alpha_{m-1}(0) = \beta_0(0) = 0, 
\alpha_m(0) = \beta_1(0) = 1$. 
Then by a family of radial transformations, the family is transformed to 
\begin{align*}
u(t,\lambda) & =  \alpha'_0(\lambda) + \alpha'_2(\lambda)t^2 + 
\cdots + \alpha'_{m-1}(\lambda)t^{m-1} + t^m, \\
w(t,\lambda) & =  \beta'_0(\lambda) + t, 
\end{align*}
for some functions 
$\alpha'_0(\lambda), \alpha'_2(\lambda), \dots, 
\alpha'_{m-1}(\lambda), \beta'_0(\lambda)$.  
The latter curve is expressed as 
$$
u = t_0(\lambda) + t_1(\lambda)w + \cdots + t_{m-1}(\lambda)w^{m-1} 
+ w^m. 
$$
By a family of linear transformations
$(u, w) \mapsto (u - t_1(\lambda)w, w)$,  it is reduced to
$$
u = t_0(\lambda) + t_2(\lambda)w^2 + \cdots + t_{m-1}(\lambda)w^{m-1} + w^m. 
$$
\enr

\ 

In general, any parametrized curve $(u(t), w(t))$ through the origin in $\C$ 
is equivalent by radial transformations 
and re-parametrizations to 
$$
u(t) = t^m + O(t^{m+1}), \quad w(t) = t^n, 
$$
for some integers $m, n$ with $m > n$. 
We have $n = 1$ for regular curves. If 
$n \geq 2$, then we call the curve an 
$(n, m)$-{\it cusp}. 
A $(2, 3)$-cusp is called a {\it simple cusp}\/ briefly a {\it cusp}. 

\bel
{\rm (The radial classification of simple cusps): } 
\label{simple cusp} 
Any simple cusp is equivalent by radial transformations 
and re-parametrizations to 
$$
u(t) = t^3, \quad w(t) = t^2. 
$$
Thus, any two simple cusps are radially equivalent to each other. 
\enl

\ 

Lemma~\ref{simple cusp} implies Proposition~\ref{Whitney cusp} 
since the discriminant of $\psi(x,y,z)=x^3+xy+z^2+\im y$ is a simple cusp. 

\ 

\noindent{\it 
Proof of Lemma~\ref{simple cusp}}\/ : 
Let $u(t) = t^3 + O(t^4), w(t) = t^2$ be a simple cusp. 
By a linear transformation on the $(u, w)$-plane and a re-parametrization of $t$, we may suppose 
the curve is given by 
$u(t) = t^3$, $w(t) = t^2 + O(t^4)$. 
Set $w(t) = t^2a(t)$ for a smooth function $a(t)$. 
Then $a(0) = 1, a'(0) = 0$. Then there exists
a smooth function $\rho(x, y)$ such that $a(t) = \rho(t^2, t^3)$ 
by the preparation theorem (\cite{GG}). 
Then the curve is radially equivalent to the curve 
$$
u(t) = \dfrac{1}{\rho(t^2, t^3)^3}t^3, \quad 
w(t) = \dfrac{1}{\rho(t^2, t^3)^2}t^2, 
$$
which is radially equivalent to $u(t) = t^2, w(t) = t^3$. 
\QED

\ber \label{3-web}
{\rm 
 Let $C_{a,b}(t)=\left(u(t,a,b),w(t,a,b)\right)$ 
be a generic two parameter family of simple cusps. 
 For each $(a,b)$, we draw tangent lines to the simple cusp $C_{a,b}$ 
from the origin. 
Then there exists a non-void open subset $U$ 
such that for $(a,b)\in U$, there are exactly three tangent rays. 
 By the assignment of the corresponding tangent points, 
we have three functions $\lambda_1,\ \lambda_2,\ \lambda_3$ on $U$; 
$C_{a,b}(\lambda_1),\ C_{a,b}(\lambda_2), C_{a,b}(\lambda_3)$ are tangent points. 
 Thus we have a triple of foliations: 
\begin{equation*}
  \lambda_1=\text{const.},\ \lambda_2=\text{const.},\ \lambda_3=\text{const.},\ 
\end{equation*}
that is, a $3$-web on $U$. 
 Moreover, radially equivalent families of simple cusps 
have isomorphic $3$-webs. 
 It is known that the classification of $3$-webs has function moduli in general
(\cite{Dufour}). 
}
\enr



%




\section{Helmholtz jet space and transversality.} \label{sec:HelmholtzJet}

We introduce the notion of Helmholtz jet spaces and show 
the transversality theorem in a Helmholtz jet space, 
as one of the main ideas to show the results in this paper. 
Note that, in \cite{IM}, analogous jet spaces are considered for 
other kinds of Monge-Amp{\` e}re equations. 

Consider the Taylor expansion around $(x, y) = (x_0, y_0)$ 
of a complex valued function $\psi$ 
on the $(x,y)$-plane: 
$$
\psi(x, y)  =  a + bX + cY + \frac{e}{2}X^2 
+ f XY + \frac{g}{2}Y^2 
+ \frac{h}{6}X^3 + \frac{k}{2}X^2Y 
+ \frac{\ell}{2}XY^2 + 
\frac{m}{6}Y^3 + \cdots. 
$$
Here we set $X = x - x_0$, $Y = y - y_0$, and $a, b, c, \dots$ are 
complex numbers. 

Suppose $\psi$ is a Helmholtz function for $k = 1$, that is, 
$\psi$ satisfies Helmholtz equation 
$\psi_{xx} + \psi_{yy} + \psi = 0$. 
Then we have 
$$
e + g + a = 0, \ 
h + \ell + b = 0, \ 
k + m + c = 0. 
$$
Therefore, we have 
$$
\begin{array}{ccl}
\vspace{0.1truecm}
\psi(x, y) & = & a + bX + cY + 
\dfrac{e}{2}X^2 + f XY - \dfrac{1}{2}(a + e)Y^2 
\\ 
 & & \quad\quad
 + \dfrac{h}{6}X^3 + \dfrac{k}{2}X^2Y - 
 \dfrac{1}{2}(b + h)XY^2 - 
\dfrac{1}{6}(c + k)Y^3 + \cdots. 
\end{array}
$$

The Taylor expansion of a function $\psi$ 
up to order $r$ around a point $(x_0, y_0)$ 
of $\R^2$ is called the $r$-{\it jet\/} of $\psi$ at $(x_0, y_0)$ 
and denoted by $j^r\psi(x_0, y_0)$. 
Denote by $J^r(\R^2, \C)$ the space of $r$-jets 
of complex valued functions on $\R^2$. 
In it, we denote by $J^r_{\rm{Helm}}(\R^2, \C)$ 
the set of $r$-jets of planar Helmholtz functions for $k = 1$: 
$$
J^r_{\rm{Helm}}(\R^2, \C) = 
\{ j^r\psi(x_0, y_0) 
\mid \psi_{xx} + \psi_{yy} + \psi = 0 {\rm{\ around \ }} 
(x_0, y_0) \}. 
$$
We call it the \emph{Helmholtz $r$-jet space}. 
For example, $J^3_{\rm{Helm}}(\R^2, \C)$ is identified with 
$\R^{16} = \R^2 \times \C \times \C^6$ with coordinates 
$x_0$, $y_0$; $a = a_1 + \im a_2$; $b = b_1 + \im b_2$, $c = c_1 + \im c_2$, 
$e = e_1 + \im e_2$, $f = f_1 + \im f_2$, $h = h_1 + \im h_2$, 
and $k = k_1 + \im k_2$. 

In general, the Taylor expansion of a Helmholtz function $\psi(x, y)$ defined 
around $(x_0, y_0)$ 
is determined by $\psi(x, y_0)$ and $\psi_y(x, y_0)$. 
Moreover,  for any given complex valued 
analytic functions $\psi_0(x)$ and 
$\psi_1(x)$ defined around $x_0$, there exists uniquely 
a complex valued function 
$\psi(x, y)$ defined around $(x_0, y_0)$ 
satisfying the Helmholtz equation, $\psi(x, y_0) = \psi_0(x)$ and 
$\psi_y(x, y_0) = \psi_1(x)$. 
The $r$-jet of $\psi$ at $(x_0, y_0)$ is determined by 
the $r$-jet of $\psi_0(x)$ at $x_0$ 
and $(r-1)$-jet of $\psi_1(x)$ at $x_0$. 
Thus, $J^r_{\rm{Helm}}(\R^2, \C)$ is identified with $\R^N$ 
for some natural number $N$. 
With any Helmholtz function $\psi$ defined around $(x_0, y_0)$, 
there is associated 
a mapping 
$$
j^r\psi \colon (\R^2, (x_0, y_0)) \to J^r_{\rm{Helm}}(\R^2, \C)
$$ 
defined by taking the $r$-jet of $\psi$ at $(x, y)$ 
for each $(x,y)$ near $(x_0, y_0)$. 
It is called the \emph{$r$-jet extension}\/ of $\psi$. 
Moreover, to any family 
$\Psi(x, y, \lambda)\colon \R^2\times\R^{\ell} \to \C$ of Helmholtz functions, 
there corresponds a mapping 
$$
j^r\Psi \colon (\R^2\times\R^{\ell}, (x_0, y_0, \lambda_0)) 
\to J^r_{\rm{Helm}}(\R^2, \C)
$$ 
by taking the $r$-jet of $\Psi(x, y, \lambda')$ at $(x, y)$ 
for each $(x,y)$ near $(x_0, y_0)$ and parameter $\lambda'$ near $\lambda_0$. 

By a similar proof to that of ordinary transversality theorem(\cite{GG}), 
we have: 

\bel \label{lem:Htransv}
Suppose a finite number of submanifolds $W_1, W_2, \dots$ of 
Helmholtz $r$-jet space $J^r_{\rm{Helm}}(\R^2, \C)$ are given. 
Then, any Helmholtz function $\psi(x, y)$ 
defined around $(x_0, y_0)$ 
is approximated (in $C^\infty$ topology) 
by a Helmholtz function $\tilde{\psi}(x, y)$ defined around 
$(x_0, y_0)$ such that 
the $r$-jet extension of $\tilde{\psi}(x, y)$ is transversal to 
any $W_i$. 
Moreover, any Helmholtz wave $\Psi(x, y, t)$ 
is approximated by a Helmholtz wave $\widetilde{\Psi}(x, y, t)$ 
such that $j^r\Psi$ is transversal to any $W_i$. 
\enl

By using Lemma~\ref{lem:Htransv}, we show 
Theorem~\ref{thm:Helmholtz_generic_fct} and Theorem~\ref{Helmholtz generic}. 
In $J^3_{\rm{Helm}}(\R^2, \C)$ with coordinates 
$x_0$, $y_0$; $a = a_1 + \im a_2$; $b = b_1 + \im b_2$, $c = c_1 + \im c_2$, 
$e = e_1 + \im e_2$, $f = f_1 + \im f_2$, $h = h_1 + \im h_2$ 
and $k = k_1 + \im k_2$, we set 
$$
W_1^6 = \{ a = 0,\ b = 0,\ c = 0 \}, 
$$
which is of codimension $6$. 
Moreover, we define a submanifold $W_3^4$ of codimension $4$ by the equations 
$$
\left\vert
\begin{array}{cc}
\vspace{0.4truecm}
{
\left\vert
\begin{array}{cc}
e_1 & c_1 \\
e_2 & c_2 
\end{array}
\right\vert
+ 
\left\vert
\begin{array}{cc}
b_1 & f_1 \\
b_2 & f_2 
\end{array}
\right\vert
}
\ 
& 
\ 
{
\left\vert
\begin{array}{cc}
f_1 & c_1 \\
f_2 & c_2 
\end{array}
\right\vert
- 
\left\vert
\begin{array}{cc}
b_1 & a_1 + e_1 \\
b_2 & a_2 + e_2 
\end{array}
\right\vert
}
\\
b_1 & c_1 
\end{array}
\right\vert = 0, 
$$
and 
$$
\left\vert
\begin{array}{cc}
\vspace{0.4truecm}
{
\left\vert
\begin{array}{cc}
e_1 & c_1 \\
e_2 & c_2 
\end{array}
\right\vert
+ 
\left\vert
\begin{array}{cc}
b_1 & f_1 \\
b_2 & f_2 
\end{array}
\right\vert
}
\ 
& 
\ 
{
\left\vert
\begin{array}{cc}
f_1 & c_1 \\
f_2 & c_2 
\end{array}
\right\vert
- 
\left\vert
\begin{array}{cc}
b_1 & a_1 + e_1 \\
b_2 & a_2 + e_2 
\end{array}
\right\vert
}
\\
b_2 & c_2 
\end{array}
\right\vert = 0, 
$$
together with $a = 0$, $b_1c_2 - b_2c_1 = 0$, 
minus a locus $W_2^5$ of more degenerate singularities, 
which is of codimension $\geq 5$. 
The definition of $W_3^4$ is from the idea of the iterated Jacobian \cite{FI}. 
Further, we set 
$$
W_4^3 = \{ a = 0,\ b_1c_2 - b_2c_1 = 0 \} \setminus (W_1^6 \cup W_2^5 \cup W_3^4), 
$$
which is of codimension $3$, 
and
$$
W_5^2 = \{ a = 0 \} \setminus (W_1^6 \cup W_2^5 \cup W_3^4 \cup W_4^3), 
$$
which is of codimension $2$. 

From Lemma~\ref{lem:Htransv}, any Helmholtz function $\psi(x, y)$ 
is approximated to a Helmholtz function 
whose $r$-jet extension is transversal to the above submanifolds 
in $J^3_{\rm{Helm}}(\R^2, \C)$. 
This implies the following from the Whitney theory 
(\cite{GG}). 
The transversality of $j^3\psi$ at $(x_0, y_0)$ to $W_5^2$ 
implies that $\psi$ has the regular phase singularity at $(x_0, y_0)$. 
Similarly, the transversality to $W_4^3$ implies the fold singularity, 
and $W_3^4$ the cusp singularity as a mapping from a plane to a plane. 
When a function has fold singularity, generically, there are two possibilities 
of phase singularities: hyperbolic and elliptic singularities. 
However, from Proposition~\ref{prop:realize_ellip}, 
there is no elliptic phase singularity for Helmholtz functions. 
This shows Theorem~\ref{thm:Helmholtz_generic_fct}. 
Furthermore, Lemma~\ref{lem:Htransv} claims 
that the transversality theorem holds even for Helmholtz waves. 
Therefore, Theorem~\ref{Helmholtz generic} is proved in the same way as above. 

\ber
The transversality to $W_2^5$ and $W_1^6$ 
means that $j^3\psi$ does not intersect to $W_2^5$ and $W_1^6$, 
for the two parameter family of Helmholtz functions. 
\enr

\section{Bifurcation problem of phase singularities of non-linear 
Schr{\" o}dinger waves. } \label{sec:Schrodinger}

We can apply our method to study phase singularities appearing in
non-linear waves. Actually we treat local analytic waves or formal  
waves.
Note that we would have to find other methods
for the study of global structure of phase singularities of non-linear waves
due to the existence of soliton solutions (see for instance \cite{Lamb}),


Let 
$$
\im\Psi_t + \dfrac{1}{2}\Psi_{xx} + f(x, \Psi) = 0 
$$
be a Schr{\" o}dinger equation for a complex valued function 
$\Psi = \Psi(t, x)$. Here $f$ is a real analytic function 
on $\R\times \C$, for instance, $f(x, \Psi) = \vert\Psi\vert^2\Psi$. 

Set $\Psi = u + \mathrm{i}w$. Then, in the case 
$f(x,\Psi) = \vert\Psi\vert^2\Psi$ 
the equation reads 
$$
\left\{
\begin{array}{ccc}
u_{xx} & = & 2w_t - 2(u^2 + w^2)u 
\vspace{0.2truecm}
\\
w_{xx} & = & - 2u_t - 2(u^2 + w^2)w
\end{array}
\right. 
$$

Let us denote by 
$$
J^r_S(\R^2, \C) := 
\{ j^r\Psi(t_0, x_0) \mid \mathrm{i}\Psi_t + \dfrac{1}{2}\Psi_{xx} 
+ f(x, \Psi) = 0 \}
$$
the {\it Schr{\" o}dinger jet space}. Then we see $J^r_S(\R^2, \C)$ 
is a submanifold of 
$J^r(\R^2, \C)$. 

In $J^r_S(\R^2, \C)$, 
the condition $\Psi = 0$ gives a smooth submanifold of $J^r_S(\R^2, \C)$ 
of codimension $2$. 
Therefore, the phase singularities appear on the $(t, x)$-plane 
at isolated points and, on  
the $x$-line, the phase singularities appear momentarily. 
The fold singularities form a submanifold 
of codimension $3$ in $J^2_S(\R^2, \C)$. The fold locus is tangent 
to the $x$-line in codimension $4$. 

The condition $\Psi = \Psi_x = 0$ also gives a smooth submanifold 
of $J^r_S(\R^2, \C)$ of codimension $4$. 
Therefore, degenerate phase singularities appear in a generic 
two parameter family of solutions, where 
$\Psi = \Psi_x = 0$. 
Consider the condition $\Psi = \Psi_x = \Psi_{xx} = 0$. If 
$\Psi(t_0, x_0) = \Psi_x(t_0, u_0) = 0$, 
then the condition 
$\Psi_{xx}(t_0, x_0) = 0$ is equivalent to that 
$\mathrm{i}\Psi_t(t_0, x_0) + f(x_0, 0) = 0$. If $f(x_0, 0) = 0$, for instance 
if $f(x,\Psi) = \vert\Psi\vert^2\Psi$, then 
$\Psi_{xx}(t_0, x_0) = 0$ if and only if $\Psi_t(t_0, x_0) = 0$. 
Thus any bifurcation on $t$ of phase singularity with $\Psi(t_0, x_0) = 
\Psi_x(t_0, x_0) = \Psi_{xx}(t_0, x_0) = 0$ is degenerate ($\Psi_t(t_0, x_0) = 0$).

\

\begin{flushleft}

Jiro ADACHI \\
Department of Mathematics, Hokkaido University, 
Sapporo 060-0810, Japan.
\vspace{-0.25truecm}
\begin{verbatim}
E-mail : j-adachi@math.sci.hokudai.ac.jp
\end{verbatim} 

\smallskip

Go-o ISHIKAWA \\ 
Department of Mathematics, Hokkaido University, 
Sapporo 060-0810, Japan.
\vspace{-0.25truecm}
\begin{verbatim}
E-mail : ishikawa@math.sci.hokudai.ac.jp
E-mail : ishikawa@topology.coe.hokudai.ac.jp
\end{verbatim} 
\end{flushleft}

\end{document}